\documentclass[12pt]{article}

\usepackage{amsmath}
\usepackage{amssymb,latexsym,xcolor,dsfont,enumitem,mathrsfs,stmaryrd,bm,mathtools} 
\usepackage{cite}
\usepackage[colorlinks=true,linkcolor=blue!50!black,pagecolor=blue!50!black,citecolor=green!50!black]{hyperref}
\usepackage[auth-sc]{authblk}
\usepackage{multirow}
\usepackage{longtable}

\usepackage{tikz}
\usepackage{pgf}
\usetikzlibrary{calc}
\usetikzlibrary{intersections}
\usetikzlibrary{fadings}

\usepackage[cp1251]{inputenc}
\usepackage[T1]{fontenc}
\usepackage[russian,english]{babel}

\makeatletter
\def\thmstyle{\it} 
\def\@begintheorem#1#2{\it \trivlist \item[\hskip
        \labelsep{\bf #1\ #2.}]\thmstyle}
\def\@opargbegintheorem#1#2#3{\it \trivlist \item[\hskip
        \labelsep{\bf #1\ #2\ (#3).}]\thmstyle}
\makeatother

\newtheorem{theorem}{{\indent}Theorem} 

\newcommand{\implication}               {\to}
  
\newcommand{\conjunction}             {\wedge}   

\makeatletter
\newcommand{\mref}{\@ifnextchar({\mref@i}{\mref@i({\Box},{p})}}
\def\mref@i(#1,#2){#1#2 \implication #2}
\makeatother
    \makeatletter
    \newcommand{\mrefp}{\@ifnextchar({\mrefp@i}{\mrefp@i({p})}}
    \def\mrefp@i(#1){\mref(\Box,#1)}
    \makeatother
\makeatletter
\newcommand{\FOref}{\@ifnextchar({\FOref@i}{\FOref@i({x},{P})}}
\def\FOref@i(#1,#2){\forall #1\,#2(#1,#1)}
\makeatother
    \makeatletter
    \newcommand{\FOrefp}{\@ifnextchar({\FOrefp@i}{\FOrefp@i({P})}}
    \def\FOrefp@i(#1){\FOref(x,#1)}
    \makeatother
\makeatletter
\newcommand{\FOrefi}{\@ifnextchar({\FOrefi@i}{\FOrefi@i({x},{P})}}
\def\FOrefi@i(#1,#2){\forall #1\,#1#2#1}
\makeatother
    \makeatletter
    \newcommand{\FOrefip}{\@ifnextchar({\FOrefip@i}{\FOrefip@i({P})}}
    \def\FOrefip@i(#1){\FOrefi(x,#1)}
    \makeatother

\makeatletter
\newcommand{\mtra}{\@ifnextchar({\mtra@i}{\mtra@i({\Box},{p})}}
\def\mtra@i(#1,#2){#1#2 \implication #1#1#2}
\makeatother
    \makeatletter
    \newcommand{\mtrap}{\@ifnextchar({\mtrap@i}{\mtrap@i({p})}}
    \def\mtrap@i(#1){\mtra(\Box,#1)}
    \makeatother
\makeatletter
\newcommand{\FOtra}{\@ifnextchar({\FOtra@i}{\FOtra@i({x},{y},{z},{P})}}
\def\FOtra@i(#1,#2,#3,#4){\forall #1\forall #2\forall #3\,(#4(#1,#2)\conjunction #4(#2,#3) \implication #4(#1,#3))}
\makeatother
    \makeatletter
    \newcommand{\FOtrap}{\@ifnextchar({\FOtrap@i}{\FOtrap@i({P})}}
    \def\FOtrap@i(#1){\FOtra(x,y,z,#1)}
    \makeatother
\makeatletter
\newcommand{\FOtrai}{\@ifnextchar({\FOtrai@i}{\FOtrai@i({x},{y},{z},{P})}}
\def\FOtrai@i(#1,#2,#3,#4){\forall #1\forall #2\forall #3\,(#1#4#2\conjunction #2#4#3 \implication #1#4#3)}
\makeatother
    \makeatletter
    \newcommand{\FOtraip}{\@ifnextchar({\FOtraip@i}{\FOtraip@i({P})}}
    \def\FOtraip@i(#1){\FOtrai(x,y,z,#1)}
    \makeatother

\makeatletter
\newcommand{\msym}{\@ifnextchar({\msym@i}{\msym@i({\Box},{\Diamond},{p})}}
\def\msym@i(#1,#2,#3){#3 \implication #1#2#3}
\makeatother
    \makeatletter
    \newcommand{\msymp}{\@ifnextchar({\msymp@i}{\msymp@i({p})}}
    \def\msymp@i(#1){\msym(\Box,\Diamond,#1)}
    \makeatother
\makeatletter
\newcommand{\FOsym}{\@ifnextchar({\FOsym@i}{\FOsym@i({x},{y},{P})}}
\def\FOsym@i(#1,#2,#3){\forall #1\forall #2\,(#3(#1,#2)\implication #3(#2,#1))}
\makeatother
    \makeatletter
    \newcommand{\FOsymp}{\@ifnextchar({\FOsymp@i}{\FOsymp@i({P})}}
    \def\FOsymp@i(#1){\FOsym(x,y,#1)}
    \makeatother
\makeatletter
\newcommand{\FOsymi}{\@ifnextchar({\FOsymi@i}{\FOsymi@i({x},{y},{P})}}
\def\FOsymi@i(#1,#2,#3){\forall #1\forall #2\,(#1#3#2 \implication #2#3#1)}
\makeatother
    \makeatletter
    \newcommand{\FOsymip}{\@ifnextchar({\FOsymip@i}{\FOsymip@i({P})}}
    \def\FOsymip@i(#1){\FOsymi(x,y,#1)}
    \makeatother

\makeatletter
\newcommand{\mser}{\@ifnextchar({\mser@i}{\mser@i({\Box},{\Diamond},{p})}}
\def\mser@i(#1,#2,#3){#1#3 \implication #2#3}
\makeatother
    \makeatletter
    \newcommand{\mserp}{\@ifnextchar({\mserp@i}{\mserp@i({p})}}
    \def\mserp@i(#1){\mser(\Box,\Diamond,#1)}
    \makeatother
\makeatletter
\newcommand{\FOser}{\@ifnextchar({\FOser@i}{\FOser@i({x},{y},{P})}}
\def\FOser@i(#1,#2,#3){\forall #1\exists #2\,#3(#1,#2)}
\makeatother
    \makeatletter
    \newcommand{\FOserp}{\@ifnextchar({\FOserp@i}{\FOserp@i({P})}}
    \def\FOserp@i(#1){\FOser(x,y,#1)}
    \makeatother
\makeatletter
\newcommand{\FOseri}{\@ifnextchar({\FOseri@i}{\FOseri@i({x},{y},{P})}}
\def\FOseri@i(#1,#2,#3){\forall #1\exists #2\,#1#3#2}
\makeatother
    \makeatletter
    \newcommand{\FOserip}{\@ifnextchar({\FOserip@i}{\FOserip@i({P})}}
    \def\FOserip@i(#1){\FOseri(x,y,#1)}
    \makeatother

\makeatletter
\newcommand{\mla}{\@ifnextchar({\mla@i}{\mla@i({\Box},{p})}}
\def\mla@i(#1,#2){#1(#1#2 \implication #2) \implication #1#2}
\makeatother
    \makeatletter
    \newcommand{\mlap}{\@ifnextchar({\mlap@i}{\mlap@i({p})}}
    \def\mlap@i(#1){\mla(\Box,#1)}
    \makeatother

\makeatletter
\newcommand{\mgrz}{\@ifnextchar({\mgrz@i}{\mgrz@i({\Box},{p})}}
\def\mgrz@i(#1,#2){#1(#1(#2 \implication #1#2) \implication #2) \implication #2}
\makeatother
    \makeatletter
    \newcommand{\mgrzp}{\@ifnextchar({\mgrzp@i}{\mgrzp@i({p})}}
    \def\mgrzp@i(#1){\mgrz(\Box,#1)}
    \makeatother

\makeatletter
\newcommand{\mwgrz}{\@ifnextchar({\mwgrz@i}{\mwgrz@i({\Box},{p})}}
\def\mwgrz@i(#1,#2){#1^+(#1(#2 \implication #1#2) \implication #2) \implication #2}
\makeatother
    \makeatletter
    \newcommand{\mwgrzp}{\@ifnextchar({\mwgrzp@i}{\mwgrzp@i({p})}}
    \def\mwgrzp@i(#1){\mwgrz(\Box,#1)}
    \makeatother

\begin{document}

\title{Logics with the axiom of convergence:\protect\\ complexity with a small number of variables\protect\\ in the language\protect\\
(extended version)\footnote{The paper was prepared by participants of the Internship Program for Employees and Postgraduate Students of Russian Educational and Scientific Institutions at the HSE University, Faculty of Mathematics, based on data obtained during the internship period. It~is an extended version of the abstracts submitted to the conference~Smirnov Readings~\cite{RybakovShcherbakov:2025:SR}.}\\~}

\author[1]{M.\,Rybakov}
\author[2]{M.\,Shcherbakov}
\affil[1]{HSM MIPT, HSE University, Tver State University}
\affil[2]{Saint-Petersburg State University, HSE University}

\date{}

\newcommand{\logic}[1]{\mathbf{#1}}
\newcommand{\num}[1]{\mathds{#1}}
\newcommand{\numbers}         [1]       {\mathds{#1}}
\newcommand{\numN}                      {\numbers{N}}
\newcommand{\numNp}                     {\numbers{N}^+}
\newcommand{\numQ}                      {\numbers{Q}}
\newcommand{\numR}                      {\numbers{R}}
\newcommand{\numC}                      {\numbers{C}}
\newcommand{\numA}                      {\numbers{A}}
\newcommand{\numAbar}                   {\bar{\numbers{A}}}
\newcommand{\ccls}[1]{\mathrm{#1}}
\newcommand{\HideText}[1]{}

\maketitle


\begin{abstract}
It is known that many modal and superintuitionistic logics are $\ccls{PSPACE}$-hard in languages with a small number of variables; however, questions about the complexity of similar fragments of many logics obtained by adding various axioms to ``standard'' ones remain unexplored. We investigate the complexity of fragments of modal logics obtained by adding an axiom requiring the convergence of the accessibility relation in Kripke frames: $\logic{S4.2}$, $\logic{K4.2}$, $\logic{Grz.2}$, and $\logic{GL.2}$. The main result is that $\logic{S4.2}$ and $\logic{Grz.2}$ are $\ccls{PSPACE}$-complete in a language with two variables, while $\logic{K4.2}$ and $\logic{GL.2}^\ast$ (a logic near to $\logic{GL.2}$) are $\ccls{PSPACE}$-complete in a language with one variable. The obtained results are extended to infinite classes of logics.
\end{abstract}

Our goal is to describe the algorithmic properties of fragments of the logics $\logic{S4.2}$, $\logic{K4.2}$, $\logic{Grz.2}$ (see, e.g.,~\cite{Chagrov:1997}), $\logic{K4.2}^\ast$ and $\logic{GL.2}^\ast$ (defined below) in languages with a small number of variables.\footnote{In the first version of the article, it was claimed that the result presented here for $\logic{GL.2}^\ast$ holds for $\logic{GL.2}$, which is not the case. We express our sincere gratitude to Taishi Kurahashi for identifying this error and for bringing to our attention that $\logic{GL.2}=\logic{GL.3}\oplus\Box\Box\bot$.} Semantically, each of these logics is characterized by some class of Kripke frames where the accessibility relation is convergent, but to be more precise, it is directed.
Here, some clarifications are needed.
In first-order language, i.e., over the logic $\logic{QCl}$, there are similar conditions:
\[
\begin{array}{lcl}
\bm{conv} & = & \forall x,y,z\,(xRy \wedge xRz \phantom{{}\wedge y\ne z} \to \exists u\, (yRu \wedge zRu)); \\
\bm{dir}  & = & \forall x,y,z\,(xRy \wedge xRz\wedge y\ne z \to \exists u\, (yRu \wedge zRu)).
\end{array}
\]
These are called the conditions of convergence and directedness, respectively.
We will actually be interested in the directedness condition, which coincides with the convergence condition in the case of reflexivity.
This condition is illustrated graphically in the figure below, where filled circles correspond to irreflexive elements, and unfilled ones to reflexive elements:
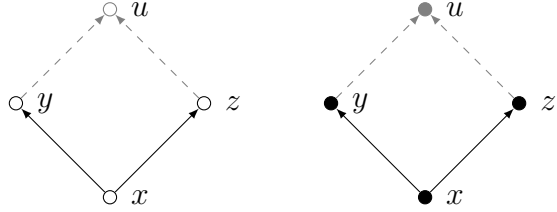
\begin{center}
\begin{tikzpicture}[scale=1.25,
   node distance=5mm,
   terminal/.style={
   rectangle,minimum size=6mm,rounded corners=3mm,
   draw=black!50,
   }
   ]

\coordinate (x)    at (  0, 0);
\coordinate (y)    at ( -1, 1);
\coordinate (z)    at ( +1, 1);
\coordinate (u)    at (  0, 2);

\draw [] (x) circle [radius=2.0pt] ;
\draw [] (y) circle [radius=2.0pt] ;
\draw [] (z) circle [radius=2.0pt] ;
\draw [color = black!50] (u) circle [radius=2.0pt] ;

\begin{scope}[>=latex]
\draw [->, shorten >= 2.75pt, shorten <= 2.75pt, dashed, color = black!50] (y)  -- (u);
\draw [->, shorten >= 2.75pt, shorten <= 2.75pt, dashed, color = black!50] (z)  -- (u);
\draw [->, shorten >= 2.75pt, shorten <= 2.75pt] (x) -- (y);
\draw [->, shorten >= 2.75pt, shorten <= 2.75pt] (x) -- (z);
\end{scope}

\node [right] at (x) {~$x$};
\node [right] at (y) {~$y$};
\node [right] at (z) {~$z$};
\node [right] at (u) {~$u$};

\end{tikzpicture}
~~~~~
\begin{tikzpicture}[scale=1.25,
   node distance=5mm,
   terminal/.style={
   rectangle,minimum size=6mm,rounded corners=3mm,
   draw=black!50,
   }
   ]

\coordinate (x)    at (  0, 0);
\coordinate (y)    at ( -1, 1);
\coordinate (z)    at ( +1, 1);
\coordinate (u)    at (  0, 2);

\draw [fill] (x) circle [radius=2.0pt] ;
\draw [fill] (y) circle [radius=2.0pt] ;
\draw [fill] (z) circle [radius=2.0pt] ;
\draw [fill, color = black!50] (u) circle [radius=2.0pt] ;

\begin{scope}[>=latex]
\draw [->, shorten >= 2.75pt, shorten <= 2.75pt, dashed, color = black!50] (y)  -- (u);
\draw [->, shorten >= 2.75pt, shorten <= 2.75pt, dashed, color = black!50] (z)  -- (u);
\draw [->, shorten >= 2.75pt, shorten <= 2.75pt] (x) -- (y);
\draw [->, shorten >= 2.75pt, shorten <= 2.75pt] (x) -- (z);
\end{scope}

\node [right] at (x) {~$x$};
\node [right] at (y) {~$y$};
\node [right] at (z) {~$z$};
\node [right] at (u) {~$u$};


\end{tikzpicture}
\end{center}
Over intuitionistic logic $\logic{Int}$, the condition of convergence (or directedness) is described by the weak law of excluded middle: 
\[
\begin{array}{lcl}
\bm{wem} & = & \neg p\vee\neg\neg p.
\end{array}
\]
Over modal logics $\logic{S4}$ and $\logic{Grz}$, the condition of convergence (or directedness) is described by the Geach axiom
\[
\begin{array}{lcl}
\bm{ga} & = & \Diamond\Box p\to\Box\Diamond p.
\end{array}
\]
Over modal logics $\logic{K4}$ and $\logic{GL}$, the directedness condition differs from the convergence condition and is described by the formula
\[
\begin{array}{lcl}
\bm{dir} & = & \Diamond(p\wedge \Box q)\to\Box(p\vee \Diamond q).
\end{array}
\]
Also, let
\[
\begin{array}{lcl}
\bm{dir}^\ast & = & \Diamond(p\wedge \Box^+ q)\to\Box(p\vee \Diamond^+ q),~~~~~~
\end{array}
\]
where $\Box^+\varphi=\varphi\wedge\Box\varphi$ and $\Diamond^+\varphi=\neg\Box^+\neg\varphi$.

Now let us define the logics that will be important to us.
Let $\logic{K}$ be the minimal normal modal propositional logic obtained by closing the set $\logic{Cl}\cup\{\Box(p\to q)\to (\Box p\to\Box q)\}$ under $\mathit{MP}$, $\mathit{Subst}$, and $\mathit{Nec}$, where $\logic{Cl}$ is the classical propositional logic:
\[
\begin{array}{lcl}
\logic{K} & = & \logic{Cl}\oplus \Box(p\to q)\to (\Box p\to\Box q).\!
\end{array}
\] 
Let, as usual,
$$
\begin{array}{lcrcl}
  \logic{KC}     & = & \logic{Int}& \!\!\!+\!\!\!      & \bm{wem};\\
  \logic{D}      & = & \logic{K}  & \!\!\!\oplus\!\!\! & \mser;\\
  \logic{T}      & = & \logic{K}  & \!\!\!\oplus\!\!\! & \mref;\\
  \logic{KTB}    & = & \logic{T}  & \!\!\!\oplus\!\!\! & \msym;\\
  \logic{K4}     & = & \logic{K}  & \!\!\!\oplus\!\!\! & \mtra;\\
  \logic{S4}     & = & \logic{K4} & \!\!\!\oplus\!\!\! & \mref;\\
  \logic{Grz}    & = & \logic{S4} & \!\!\!\oplus\!\!\! & \mgrz;\\
  \logic{GL}     & = & \logic{K4} & \!\!\!\oplus\!\!\! & \mla;\\
  \logic{wGrz}   & = & \logic{K4} & \!\!\!\oplus\!\!\! & \mwgrz;\\
  \logic{K4.2}   & = & \logic{K4}   & \!\!\!\oplus\!\!\! & \bm{dir};\\
  \logic{S4.2}   & = & \logic{S4}   & \!\!\!\oplus\!\!\! & \bm{dir} \quad =~~~~\logic{S4} ~\oplus~ \bm{ga};\\
  \logic{Grz.2}  & = & \logic{Grz}  & \!\!\!\oplus\!\!\! & \bm{dir} \quad =~~\logic{Grz} ~\oplus~ \bm{ga};\\
  \logic{GL.2}   & = & \logic{GL}   & \!\!\!\oplus\!\!\! & \bm{dir};\\
  \logic{wGrz.2} & = & \logic{wGrz} & \!\!\!\oplus\!\!\! & \bm{dir},
\end{array}
$$
and let
$$
\begin{array}{lcrcl}
  \logic{K4.2}^\ast   & = & \logic{K4}   & \!\!\!\oplus\!\!\! & \bm{dir}^\ast; 
  \phantom{\quad\! =~~~~~~\logic{K4}\hfill ~\oplus~ \bm{ga}^\ast}
  \hspace{0.5em}\\
  \logic{GL.2}^\ast   & = & \logic{GL}   & \!\!\!\oplus\!\!\! & \bm{dir}^\ast; 
  \phantom{\quad\! =~~~~~\logic{GL}\hfill ~\oplus~ \bm{ga}^\ast}
  \hspace{0.5em}\\
  \logic{wGrz.2}^\ast & = & \logic{wGrz} & \!\!\!\oplus\!\!\! & \bm{dir}^\ast. 
  \phantom{\quad\! =~~\logic{wGrz}\hfill ~\oplus~ \bm{ga}^\ast}
  \hspace{0.5em}\\
\end{array}
$$
It is not difficult to see that the following inclusions hold: $\logic{K4.2}^\ast \subset \logic{K4.2}$, $\logic{GL.2}^\ast \subset \logic{GL.2}$, and $\logic{wGrz.2}^\ast \subset \logic{wGrz.2}$. Also, $\logic{GL.2}=\logic{GL.3}\oplus\Box\Box\bot$, and hence, tabular; then, it is $\ccls{coNP}$-complete, and, for every $n\in\numN$, the $n$-variable fragment of $\logic{GL.2}$ is decidable in polynomial time.

Using methods for studying the complexity of modal and superintuitionistic logics~\cite{Lad:1977,Statman:1979,Shapirovsky:2018,Shapirovsky:2022,Chagrov:1985}, it is easy to show that the logics under consideration are $\ccls{PSPACE}$-complete in a language with a countable set of propositional variables. At the same time, it is known that for many ``standard'' logics with $\ccls{PSPACE}$-complete decision problems, their fragments with a small number of variables are $\ccls{PSPACE}$-complete as well: in the case of $\logic{Int}$ and $\logic{KC}$, two variables suffice~\cite{MR:2004,MR:2007}, in the case of $\logic{D}$, $\logic{T}$, $\logic{KTB}$, $\logic{S4}$, $\logic{GL}$, and $\logic{Grz}$ one variable suffices, and in the case of $\logic{K}$, $\logic{K4}$, and $\logic{wGrz}$ no variables are needed at all~\cite{Halpern:1995,Svejdar:2003,ChR:2003:AiML,AgadzhanianRybakov:2022:arXiv,Spaan:1993:PhD,RSh:2019:IGPL}. Note that for all these logics, further reducing the number of variables in the language (if possible) leads to polynomially decidable fragments: for example, the one-variable fragments of all superintuitionistic logics are polynomially decidable due to the construction by L.\,Rieger and I.\,Nishimura~\cite{Rieger:1952,Nishimura:1960,Riger:1957}, and in extensions of the logic $\logic{D}$, every variable-free formula is equivalent to either $\bot$ or $\top$, which gives polynomial decidability of the variable-free fragments of all such logics.

Since $\logic{S4.2}$ and $\logic{Grz.2}$ are the modal companions of the logic of the weak law of excluded middle $\logic{KC}$, it immediately follows from~\cite{MR:2004,MR:2007} that the two-variable fragments of $\logic{S4.2}$ and $\logic{Grz.2}$ are $\ccls{PSPACE}$-complete.
It is not hard to show that a similar result holds for logics $\logic{K4.2}$, $\logic{K4.2}^\ast$, $\logic{wGrz.2}$, $\logic{wGrz.2}^\ast$, and $\logic{GL.2}^\ast$. However, as noted above, to prove $\ccls{PSPACE}$-completeness of $\logic{K4}$, $\logic{wGrz}$ and $\logic{GL}$, one variable suffices (for $\logic{K4}$ and $\logic{wGrz}$ even variable-free formulas suffice); the same is true for the ``superintuitionistic'' fragments of these logics (i.e., the preimages under the G\"{o}del translation)---the basic and formal propositional logics $\logic{BPL}$ and $\logic{FPL}$ \cite{MR:2003,MR:2006}, introduced by A.\,Visser~\cite{Visser:1981}. These observations led us to conjecture that to prove $\ccls{PSPACE}$-hardness of $\logic{K4.2}$ and $\logic{GL.2}^\ast$, there is no need to use two variables, and one can suffice.

During the research, the following general results were obtained.

\begin{theorem}
Every modal logic between\/ $\logic{K}$ and\/ $\logic{Grz.2}$ is\/ $\ccls{PSPACE}$-hard in a language with two propositional variables.
\end{theorem}

\begin{theorem}
Every modal logic between\/ $\logic{K}$ and\/ $\logic{GL.2}^\ast$ or between\/ $\logic{K}$ and\/ $\logic{K4.2}$ is\/ $\ccls{PSPACE}$-hard in a language with one propositional variable.
\end{theorem}

The first theorem generalizes the observation made about the two-variable fragments of $\logic{S4.2}$ and $\logic{Grz.2}$, while the second gives $\ccls{PSPACE}$-completeness of the one-variable fragments of $\logic{K4.2}$, $\logic{K4.2}^\ast$, $\logic{wGrz.2}$, $\logic{wGrz.2}^\ast$, and $\logic{GL.2}^\ast$ (regarding $\logic{wGrz}$, see~\cite{Litak:2007}). Note that the question remains about the complexity of the variable-free fragments of $\logic{K4.2}$, $\logic{K4.2}^\ast$, $\logic{wGrz.2}$, and $\logic{wGrz.2}^\ast$, as well as the one-variable fragments of $\logic{S4.2}$ and $\logic{Grz.2}$. It would also be interesting to determine the complexity of the ``superintuitionistic'' fragments of $\logic{K4.2}$, $\logic{K4.2}^\ast$, and $\logic{GL.2}^\ast$ in languages with a finite number of variables.



\end{document}